\documentclass{gtpart}

\usepackage{graphicx}
\usepackage{pinlabel}  
\usepackage{amscd}
\usepackage{amsfonts}
\usepackage{mathrsfs}
\usepackage[all]{xy}
\usepackage{verbatim}

\newcommand{\s}{\sigma}
\newcommand{\p}[1][]{\pi_{#1}}
\newcommand{\y}{\lambda}

\newcommand{\g}{\gamma}
\renewcommand{\a}{\alpha}
\renewcommand{\p}{\rho}
\renewcommand{\t}{\tau}

\newcommand{\Diff}{\text{Diff}}
\newcommand{\Ind}{\text{Ind}}
\newcommand{\Stab}{\text{Stab}}
\newcommand{\Aut}{\text{Aut}}
\newcommand{\Inn}{\text{Inn}}
\newcommand{\rel}{\text{rel}}
\renewcommand{\dim}{\text{dim}}
\newcommand{\IA}{\text{IA}}
\newcommand{\I}{\mathbb{I}}

\newcommand{\PS}{P\Sigma}

\newcommand{\cP}{\mathscr{P}}
\newcommand{\bcP}{\overline{\mathscr{P}}}
\renewcommand{\S}{\Sigma}

\newcommand{\N}{\mathbb{N}}

\date{\today}
\title{Representation stability for the cohomology of the pure string motion groups}  
\author[J\,C\,H\,Wilson]{Jennifer C. H. Wilson}

\givenname{Jennifer}
\surname{Wilson}
\email{wilsonj@math.uchicago.edu}
\address{Department of Mathematics \nl University of Chicago \nl 5734 University Ave. \nl Chicago, IL 60637}
\urladdr{http://math.uchicago.edu/~wilsonj/}



\arxivreference{math.GT/1108.1255}
\arxivpassword{f2sq9}

\volumenumber{}
\issuenumber{}
\publicationyear{}
\papernumber{}
\startpage{}
\endpage{}
\doi{}
\MR{}
\Zbl{}
\received{}
\revised{}
\accepted{}
\published{}
\publishedonline{}
\proposed{}
\seconded{}
\corresponding{}
\editor{}
\version{}

\newtheorem{thm}{Theorem}[section]    
\newtheorem{lem}{Lemma}[section]          

\newtheorem{cor}{Corollary}[section]
\theoremstyle{definition}
\newtheorem{defn}[thm]{Definition}    
\newtheorem*{rem}{Remark}             
\makeatletter
\let\c@lem=\c@thm
\makeatother

\makeatletter
\let\c@cor=\c@thm
\makeatother

\begin{document}

\begin{abstract}

The cohomology of the pure string motion group $\PS_n$ admits a natural action by the hyperoctahedral group $W_n$. In \cite{RepStability}, Church and Farb conjectured that for each $k \geq 1$, the cohomology groups $H^k(\PS_n; \Q)$ are uniformly representation stable; that is, the description of the decomposition of $H^k(\PS_n; \Q)$ into irreducible $W_n$--representations stabilizes for $n >> k$. We use a characterization of $H^*(\PS_n; \Q)$ given by Jensen, McCammond, and Meier to prove this conjecture. Using a transfer argument, we further deduce that the rational cohomology groups of the string motion group $H^k(\S_n;\Q)$ vanish for $k \geq 1$. We also prove that the subgroup of $\S_n^+ \subseteq \S_n$ of orientation-preserving string motions, also known as the braid-permutation group, is rationally cohomologically stable in the classical sense.\\

\end{abstract}

\begin{asciiabstract}
 The cohomology of the pure string motion group PSigma_n admits a natural action by the hyperoctahedral group W_n. In recent work, Church and Farb conjectured that for each k > 0, the sequence of degree k rational cohomology groups of PSigma_n is uniformly representation stable with respect to the induced action by W_n, that is, the description of the groups' decompositions into irreducible W_n representations stabilizes for n >> k. We use a characterization of the cohomology groups given by Jensen, McCammond, and Meier to prove this conjecture. Using a transfer argument, we further deduce that the rational cohomology groups of the string motion group vanish in positive degree. We also prove that the subgroup of orientation-preserving string motions, also known as the braid-permutation group, is rationally cohomologically stable in the classical sense. 
\end{asciiabstract}

\maketitle

\section{Introduction}

 Let $C_n = \g_1 \cup \g_2 \cup \cdots \cup \g_n$ be the disjoint union of $n$ smoothly embedded, oriented, unlinked, unknotted circles $\gamma_i$ in $\R^3$. We define $\S_n$ to be the group of motions of $C_n$, as follows: A \emph{motion} of $C_n$ is a path of diffeomorphisms $f_t \in \Diff(\R^3)$ such that $f_0$ is the identity and $f_1$ stabilizes $C_n$ setwise. Two motions $f_{t,0}$ and $f_{t,1}$ are considered equivalent if they are smoothly isotopic through an isotopy $f_{t,s}$ such that $f_{0,s}$ and $f_{1,s}$ stabilize $C_n$. The product of two motions $f_t$ and $g_t$ is given by
\begin{align*}
 (g \cdot f)_t = \left\{ \begin{array}{ll}
         f_{2t} & \mbox{$0 \leq t \leq \frac{1}{2}$},\\
        g_{2t-1}\circ f_1 & \mbox{$\frac{1}{2} \leq t \leq 1$}.\end{array} \right.
\end{align*}

With this product, the set of these homotopy classes of motions forms the \emph{string motion group}, also called the \emph{circle-braid group}, with identity the \emph{stationary motion}. \\ 

The string motion group $\S_n$ is a generalization of the classical braid group $B_n$, the motions of $n$ points in the plane, though the circles defining $\S_n$ may additionally pass over and through each other. The analogue of the pure braid group $P_n$ is the \emph{pure string motion group} $\PS_n \subseteq \S_n$, the subgroup of motions in which all circles return to their original positions and orientations. \\

The groups $\S_n$ and $\PS_n$ may alternatively be defined in terms of spaces of smooth embeddings of circles into $\R^3$, or of configurations of `rings', as described in \fullref{Motions}. The structure of the group $\S_n$ is described in more detail in \fullref{PSn}, where it is identified with a subgroup of automorphisms of the free group on $n$ letters. \\

The quotient $\S_n / \PS_n$ is the \emph{hyperoctahedral group} $W_n$, the wreath product $$W_n=\Z / 2 \Z \wr S_n := (\Z / 2\Z)^n \rtimes S_n,$$ where $S_n$ conjugates $(\Z / 2\Z)^n$ by permuting the coordinates. The hyperoctahedral group is also called the \emph{signed permutation group}, since it may be described as the symmetries of the set $$S=\{ \{-1,1\}, \{-2, 2\}, \ldots, \{-n, n\} \},$$ where the $k^{th}$ factor of $(\Z/2\Z)^n$ transposes the elements in the block $\{-k, k \}$, and $S_n$ permutes the $n$ blocks. The hyperoctahedral group $W_n=\S_n / \PS_n$ acts by conjugation on the pure string motion group $\PS_n$ by relabelling and reversing the orientations of the $n$ circles. \\

The cohomology of $\PS_n$ is described in \fullref{CohomologyPSn}. It has been studied by Collins \cite{Collins}, Brownstein and Lee \cite{BrownsteinLee}, Brady, McCammond, Meier, and Miller, \cite{DualityGroup}, Pettet \cite{PettetFiniteness}, Jensen, McCammond, and Meier \cite{JMM}, and others. As with the pure braid group, the pure motion group is not cohomologically stable in the classical sense: Jensen, McCammond, and Meier determined the formula $$\dim_{\Q}\, H^k(\PS_n; \Q) = {n-1 \choose k}n^k.$$ For each fixed $k \geq 1$, the dimension grows unboundedly with $n$. \\

In \cite{RepStability}, Church and Farb define an alternate notion of stability, which they call representation stability. We may consider the group $H^k(\PS_n; \Q)$ as a rational representation of the hyperoctahedral group $W_n$ under the action of $W_n$ induced by its action on $\PS_n$. The irreducible representations of $W_n$ over $\Q$ are classified by pairs of partitions $\y = (\y^+, \y^-)$ with $|\y^+| + |\y^-| =n$, a classification described in more detail in \fullref{DefinitionRepStability}. Church and Farb conjectured that, in analogy with the pure braid group, the cohomology groups $H^k(\PS_n; \Q)$ stabilize in terms of their decompositions into irreducible representations of $W_n$  \cite[Conjecture 4.8]{RepStability}. We prove this conjecture in \fullref{MainThm}. \\

To make sense of this statement, we need a system of classifying the irreducible representations of $W_n$ that is independent of $n$. As in \cite{RepStability}, we write $$V( (\y^+_1, \y^+_2, \ldots , \y^+_p), (\y^-_0, \y^-_1, \ldots , \y^-_q) )_n$$ to denote the irreducible representation of $W_n$ associated to the double partition $(\y^+, \y^-)$ of $n$ with $$\y^+ = (n - k, \y^+_1, \y^+_2, \ldots , \y^+_p), \quad \y^-= (\y^-_0, \y^-_1, \ldots , \y^-_q ),$$ where $k= (\sum_{i=1}^p \y^+_i) + (\sum_{i=0}^q \y^-_i).$ In this notation, for example, the decomposition of $H^1(\PS_n; \Q)$ into irreducible $W_n$--representations is: 

\[
\begin{array} {rcll}
H^1(\PS_2; \Q) &=& V\big((0), (1)\big)_2 &  \\
 H^1(\PS_n; \Q) &=& V\big((0),(1)\big)_n \oplus V\big((1),(1)\big)_n & \text{for all } n \geq 3. \\
\end{array}
\]

This decomposition stabilizes once $n=3$. \\

A complete definition of representation stability is given in \fullref{DefinitionRepStability}. Using this framework, we adapt Church and Farb's proof of the uniform representation stability of the cohomology of $P_n$ to a description of $H^k(\PS_n ; \Q)$ given by Jensen, McCammond, and Meier in \cite{JMM} to prove our main result:

\newtheorem*{WnStable}{\fullref{WnStable}}
\begin{WnStable}
For each fixed $k \geq 0$, the sequence of $W_n$--representations $$\{H^k(\PS_n ; \Q)\}_{n \in \N}$$ is uniformly representation stable with respect to the maps $$\phi_n \co H^k(\PS_n ; \Q) \to H^k(\PS_{n+1} ; \Q)$$ induced by the `forgetful' map $\PS_{n+1} \to \PS_n$. The sequence stabilizes once $n \geq 4k$.
\end{WnStable}

It is generally true that the restriction of a uniformly stable sequence of $W_n$--representations to the action of $S_n \subseteq W_n$ yields a uniformly representation stable sequence of $S_n$--representations, a result that we prove in a forthcoming paper. However, by using the same method of proof as for \fullref{WnStable}, we prove uniform representation stability of the sequences $H^k(\PS_n ; \Q)$ with respect to the action of $S_n$ with an improved stable range. 

\newtheorem*{SnStable}{\fullref{SnStable}}
\begin{SnStable}
For each fixed $k \geq 0$, the sequence of groups $$\{H^k(\PS_n ; \Q)\}_{n \in \N},$$ considered as $S_n$--representations, is uniformly representation stable with respect to the maps $$\phi_n \co H^k(\PS_n ; \Q) \to H^k(\PS_{n+1} ; \Q)$$ induced by the `forgetful' map $\PS_{n+1} \to \PS_n$. The sequence stabilizes once $n \geq 4k$.
\end{SnStable}

If we define $$V( \y_1, \ldots \y_d)_n$$ to be the irreducible $S_n$--representation associated to the partition $$\y = (n - k, \y_1, \ldots, \y_d),$$ with $k= \sum_{i=1}^d \y_i$, then the decomposition of $H^1(\PS_n; \Q)$ into irreducible $S_n$--representations is: 

\[
\begin{array} {rcll}
H^1(\PS_2; \Q) &=& V(0)_2 \oplus V( 1)_2 &  \\
H^1(\PS_3; \Q) &=& V(0)_3 \oplus V( 1)_3^{\bigoplus 2} \oplus V(1,1)_3 &  \\
 H^1(\PS_n; \Q) &=& V(0)_n \oplus V(1)_n^{\bigoplus 2} \oplus V(1,1)_n \oplus V(2)_n & \text{for all } n \geq 4. \\
\end{array}
\]

The decomposition of $H^2(\PS_n; \Q)$ is:
\begin{align*}
H^2(\PS_2) =& \; 0 & \\
H^2(\PS_3) =& \; V(1,1)_3^{\bigoplus 2} \oplus V(1)_3^{\bigoplus 3}  \oplus V(0)_3 & \\
H^2(\PS_4) =& \; V(1,1,1)_4^{\bigoplus 2}  \oplus V(1,1)_4^{\bigoplus 7}  \oplus V(2)_4^{\bigoplus 3} \oplus V(1)_4^{\bigoplus 6} \oplus V(0)_4 & \\ 
H^2(\PS_5) =& \;  V(1,1,1)_5^{\bigoplus 4} \oplus V(2,1)_5^{\bigoplus 5} \oplus V(1,1)_5^{\bigoplus 9} \oplus V(2)_5^{\bigoplus 6} \oplus V(1)_5^{\bigoplus 6} \oplus V(0)_5 & \\ 
H^2(\PS_6) =& \; V(2,1,1)_6^{\bigoplus 2}  \oplus  V(1,1,1)_6^{\bigoplus 4} \oplus  V(2,1)_6^{\bigoplus 7} \oplus V(3)_6^{\bigoplus 3}  & \\ & \oplus  V(1,1)_6^{\bigoplus 9} \oplus V(2)_6^{\bigoplus 6} \oplus  V(1)_6^{\bigoplus 6} \oplus V(0)_6 & \\ 
H^2(\PS_n) =& \; V(2,1,1)_n^{\bigoplus 2} \oplus V(3,1)_n^{\bigoplus 2} \oplus V(1,1,1)_n^{\bigoplus 4} \oplus V(2,1)_n^{\bigoplus 7} \oplus V(3)_n^{\bigoplus 3} & \\ & \oplus V(1,1)_n^{\bigoplus 9} \oplus V(2)_n^{\bigoplus 6} \oplus V(1)_n^{\bigoplus 6} \oplus V(0)_n  \qquad \qquad \text{ for all }   n \geq 7.  
\end{align*}

These results provide a new perspective on the structure and growth of the cohomology groups $H^k(\PS_n;\Q)$ for large values of $n$.  \\

The structure of the $H^k(\PS_n;\Q)$ provides insight into the cohomology of the string motion group $\S_n$. A transfer argument identifies the rational cohomology $H^k(\S_n;\Q)$ with the $W_n$--invariants $H^k(\PS_n;\Q)^{W_n}$. In \fullref{StringMotionGroup}, using a combinatorial argument, we show that the trivial $W_n$--representation does not occur in $H^k(\PS_n;\Q)$ for $k>0$, and so prove:

\newtheorem*{SigmaVanish}{\fullref{SigmaVanish}}
\begin{SigmaVanish}
 For $k \geq 1$, the rational cohomology groups $H^k(\S_n; \Q)=0$. In particular, the integral homology and cohomology groups of $\S_n$ are finite. 
\end{SigmaVanish}

The subgroup $\S_n^+ \subseteq \S_n$ of string motions that preserve the orientations of the circles has been called the \emph{braid-permutation group}. The quotient $\S_n^+ / \PS_n = S_n$, and so the structure of $H^k(\PS_n;\Q)$ as an $S_n$--representation encodes the structure of $H^k(\S_n^+;\Q)$.  In \fullref{BraidPermutation}, we use a transfer argument to conclude the following corollary. 

\newtheorem*{BPStable}{\fullref{BPStable}}
\begin{BPStable}
 The braid-permutation group $\S_n^+$ of orientation-preserving string motions is rationally cohomologically stable in the classical sense. For each $k \geq 0$, $H^k(\S_n^+; \Q) \cong H^k(\S_{n+1}^+; \Q)$ once $n\geq 4k$.
\end{BPStable}

By analyzing the combinatorics of $H^k(\PS_n;\Q)$, we gain some insights into the structure of $H^k(\S_n^+; \Q)$, and we compute the stable dimensions in small degrees.

\begin{rem}
 We have recently learned that the integral (co)homology of $\S_n$ and $\S_n^+$ could be computed using the techniques of Griffin in \cite{Griffin}, which gives a general approach to computing the (co)homology of the symmetric automorphism groups of certain free products of groups. 
\end{rem}

\subsubsection*{Acknowledgments} I would like to extend my thanks to Benson Farb. His guidance and generous feedback made this work possible. I am additionally grateful to Rita Jimenez Rolland for her helpful comments on early versions of this paper, and to Allen Hatcher, for his numerous corrections and suggestions. I would like to thank Tom Church for a detailed proofreading and recommendations for improvement to the paper. \\ 

This work was funded in part by the Natural Sciences and Engineering Research Council of Canada. I am grateful for their support.

\section{Some Background on the Theory of Motion Groups} \label{Motions}

The string motion group $\S_n$ was first defined by Dahm in his 1962 thesis \cite{Dahm}, and the theory of motion groups was further developed by Goldsmith \cite{Goldsmith} and others. These groups were originally studied in the topological category, as follows: Given a manifold $M$ and a compact subspace $N$ in the interior of $M$, let $H(M)$ be the space of self-homeomorphisms of $M$ with the compact-open topology. Let $H(M,N)$ be the subspace of homeomorphisms that stabilize $N$ setwise. Let $H_C(M) \subseteq H(M)$ and $H_C(M,N) \subseteq H(M,N)$ be those homeomorphisms with compact support. Let $\I_M$ be the identity map on $M$. The \emph{motion group} $\mathcal{M}(M,N)$ is then defined as the relative fundamental group $\pi_1( H_C(M), H_C(M,N), \I_M)$. Although relative fundamental groups are not groups in general, a group structure on $\mathcal{M}(M,N)$ derives from the fact that $H_C(M,N) \subseteq H_C(M)$ is a subgroup. \\ 

These motion groups were reconceived in the differential category by Wattenberg \cite{Wattenberg}. Let $M$ be a smooth, connected manifold without boundary and $i \co N \hookrightarrow M$ the smooth embedding of a compact submanifold $N$. Let $E(N,M)$ be the space of smooth embeddings $f \co N \to M$ with the $C^{\infty}$ topology, and $P(N,M)$ the subspace of maps $f \in E(N,M)$ such that $f(N)=i(N)$. Let $D(M)$ be the space of of self-diffeomorphisms of $M$ with the $C^{\infty}$ topology, $D(M,N)$ the subspace stabilizing $i(N)$, and $D_C(M)$ and $D_C(M,N)$ the diffeomorphisms with compact support. Wattenberg proved that $$\pi_1( D(M), D(M,N), \I_M) \cong \pi_1(E(N,M), P(N,M), i).$$ When $M = \R^3$ and $N=C_n$, this is the string motion group $\S_n$. Wattenberg proved, in this case, that the group is isomorphic to $\pi_1( D_C(\R^3), D_C(\R^3,C_n), \I_{\R^3})$, and moreover coincides with Dahm's motion group $\mathcal{M}(\R^3, C_n)$.\\

In the introduction, we defined $\S_n$ as the group $\pi_1( D(\R^3), D(\R^3,C_n), \I_{\R^3})$, though by Wattenberg's results we may equivalently define $\S_n$ as the fundamental group of the space of smooth unlinks of $n$ components, that is, smooth embeddings $E(C_n, \R^3)$ relative embeddings stabilizing $C_n$. In this framework, the pure string motion group $\PS_n$ is then the fundamental group of the space of smooth, labelled, oriented unlinks of $n$ components, relative those embeddings preserving the order and orientation of each circle. We note that the hypothesis of smoothness is necessary: the analogous construction with continuous embeddings encounters problems due to the existence of `pathological' isotopies, for example, continuous paths that pull a knot through to an unknot. Such isotopies do not extend to motions of the ambient space. \\

Recently, Brendle and Hatcher have shown in \cite{BrendleHatcher} that it suffices to consider configuration spaces of `rings', that is, embeddings of Euclidean circles in $\R^3$. They define $\mathcal{PUR}(n)$ to be the configuration space of $n$ labelled, oriented, disjoint, unlinked rings. Each ring can be specified by vector determining its centre, and a `normal' vector that indicates its radius and orientation. Thus $\mathcal{PUR}(n)$ can be realized as an open subset of $\R^{6n}$. The configuration space $\mathcal{R}(n)$ of disjoint, unlinked rings is also a $6n$--manifold, the quotient of $\mathcal{PUR}(n)$ by the action of the hyperoctahedral group $W_n$. Brendle and Hatcher prove that $\S_n = \pi_1( \mathcal{R}(n) )$ and that $\PS_n = \pi_1(\mathcal{PUR}(n))$. This characterization of these groups facilitates the study of their algebraic structure. 

\section{The (Pure) String Motion Group} \label{PSn}

A valuable technique in the study of the string motion group is to embed the group in the automorphism group of a free group. \\

If we take the one-point compactification $S^3 = \R^3 \cup \{ \infty \}$ of $\R^3$, then the fundamental group $\pi_1(S^3 \setminus C_n, \infty)$ is the free group $F_n$ on $n$ generators. Let $x_i$ denote the generator linked once (in the positive sense) with $\g_i$, and unlinked with the other circles. By extending a string motion $f_t$ to $S^3$ and taking the induced map $(f_1)_*$ on the fundamental group  $\pi_1(S^3 \setminus C_n, \infty)$, we get a map $$D \co \S_n \to \Aut(F_n)$$ where $\Aut(F_n)$ denotes the automorphism group of $F_n$. The map $D$ is a particular instance of the \emph{Dahm homomorphism} on motion groups, developed by Dahm in \cite{Dahm} and described in \cite{Goldsmith} by Goldsmith. In the case of the string motion group $\S_n$, $D$ is injective, and its image is the group of \emph{symmetric automorphisms}, the automorphisms that map each generator $x_i$ to a conjugate of a generator or of the inverse of a generator. This group is generated by elements $\a_{i,j}$, $\t_i$, and $\p_i$,
 defined as follows: 
 \begin{align*}
 \a_{i,j} = \left\{ \begin{array}{ll}
        x_i \mapsto x_jx_ix_j^{-1} & \\
        x_k \mapsto x_k \quad \;  (k \neq i) &   \end{array} \right. 
\end{align*}\begin{align*}
 \t_i = \left\{ \begin{array}{ll}
        x_i \mapsto x_{i+1} & \\
        x_{i+1} \mapsto x_i & \\ 
        x_k \mapsto x_k \quad (k \neq i,\, i+1) & \end{array} \right. 
 \p_i = \left\{ \begin{array}{ll}
        x_i \mapsto x_i^{-1}& \\
        x_k \mapsto x_k \quad (k \neq i) & \end{array} \right.
\end{align*}

The automorphism $\a_{i,j}$, $i \neq j$, corresponds to a string motion where the $i^{th}$ circle $\gamma_i$ passes through the $j^{th}$ circle $\gamma_j$, and returns to its original position and orientation. The automorphism $\t_i$ is induced by a string motion that transposes the circles $\gamma_i$ with $\gamma_{i+1}$ while preserving their orientations, and $\p_i$ is induced by a string motion that reverses the orientation of the circle $\gamma_i$. The algebraic structure of this group was described by authors including McCool \cite{McCool} and Rubinzstein \cite{Rubinsztein}.\\

We note that the classical braid group $B_n$ embeds in $\S_n$ as the (infinite-index) subgroup generated by elements of the form $\a_{i,i+1} \t_i$, for $i=1, \ldots, (n-1).$ This is equivalently the group of symmetric automorphisms that fix the product $x_n x_{n-1} \cdots x_2 x_1$.\\

The pure string motion group $\PS_n$ is the group of all automorphisms that map each generator $x_i$ to a conjugate of itself, also called the \emph{basis-conjugating automorphisms} or the \emph{pure symmetric automorphisms}. This group is torsion-free, and is generated by the elements $\a_{i,j}$, $ i \neq j$. It is the intersection of $\S_n$ with the subgroup $\IA_n \subseteq \Aut(F_n)$ of automorphisms that induce the trivial map on the homology group $H^1(F_n; \Z) \cong \Z^n$. \\

The quotient $\S_n / \PS_n = \langle \p_i, \t_i \rangle$ is the hyperoctahedral group $W_n$ described in the introduction. As explained in by Brownstein and Lee in \cite{BrownsteinLee}, these groups fit in to a commutative diagram of short exact sequences, 

\begin{center}

$\xymatrix@H+15pt{ & 1 \ar[r] & \PS_n \ar@{^{(}->} [d] \ar[r] & \S_n \ar@{^{(}->} [d] \ar[r] & W_n \ar@{^{(}->} [d] \ar[r] & 1 \\ & 1 \ar[r] & \IA_n \ar[r]  & \Aut(F_n) \ar[r] & GL_n(\Z) \ar[r] & 1 } $ \\

\end{center}

The top exact sequence splits, giving a decomposition of $\S_n$ as a semidirect product $\S_n \cong \PS_n \rtimes W_n$. The group $W_n$ acts on $\PS_n$ as follows:

\[
 \p_i (\a_{j,k}) \p_i^{-1} = \left\{ \begin{array}{lll}
        \a_{j,k}^{-1} & (k=i) &  \\  
          \a_{j,k}  & (k \neq i) &  \end{array} \right. \qquad
\] \[
\t_i (\a_{j,k}) \t_i^{-1} =  \a_{ \sigma(j), \sigma(k) } \;\; \text{where } \sigma = ( i \; \; \; i+1 ) \in S_n.
\]

\section{The Cohomology of the Pure String Motion Group} \label{CohomologyPSn}

Define $[n]: = \{ 1, 2, \ldots, n \}$. McCool gave a presentation for $\PS_n$ in which all relations are commutators \cite{McCool}, and so the abelianization $H_1(\PS_n; \Z)$ of $\PS_n$ is the free abelian group $\Z[\a_{i,j}]$, with $i,j \in [n]$, $i \neq j$. The cohomology ring $H^*(\PS_n; \Z)$ is generated by the dual basis $\{ \a^*_{i,j} \}$. \\

Brownstein and Lee conjectured a presentation for the cohomology ring $H^*(\PS_n;\Z)$ \cite[Conjecture 4.6]{BrownsteinLee}, which was established by Jensen, McCammond and Meier in \cite{JMM}. Their result is: 

\begin{thm}  \cite[Theorem 6.7]{JMM}. 
 The cohomology ring $H^*(\PS_n; \Z)$ is the exterior algebra generated by the degree-one classes $\a_{i,j}^*$, with $i,j \in [n]$, $i \neq j$, modulo the relations 
\begin{align*}
      (1) \;  \a^*_{i,j} \wedge \a^*_{j,i} = 0   \qquad \qquad  (2) \; \a^*_{k,j} \wedge  \a^*_{j,i} -  \a^*_{k,j} \wedge  \a^*_{k,i} + \a^*_{i,j} \wedge  \a^*_{k,i} = 0 
\end{align*}
\end{thm}
These relations, computed in by Brownstein and Lee in \cite{BrownsteinLee}, all pull back from $H^*(\PS_3 ;\Z)$, under the map induced by the map $\PS_n \to \PS_3$ determined by `forgetting' all but the $i^{th}$, $j^{th}$, and $k^{th}$ circle. \\

The cohomology groups are computed in \cite{JMM} by Jensen, McCammond, and Meier using the Lyndon-Hochschild-Serre spectral sequence associated with the short exact sequence,

\begin{center}

$\xymatrix{ & 1 \ar[r] & \Inn(F_n) \ar[r] & \PS_n \ar[r] & \PS_n / \Inn(F_n) \ar[r] & 1} $ \\

\end{center}

 In the proof of Lemma 6.6 in \cite{JMM}, Jensen, McCammond, and Meier determine that any \emph{cyclic product} $$\a_{i,j}^* \wedge \a_{j,l}^* \wedge \a_{l,m}^*  \wedge \cdots \wedge \a_{r,s}^* \wedge \a_{s,i}^*$$ is trivial. They argue that $H^k(\PS_n;\Z)$ is generated by elements of the form $\{ \a_{i_1, j_1}^* \wedge \a_{i_2, j_2}^* \wedge \cdots \wedge \a_{i_k, j_k}^* \}$ that have distinct indices $i_m$ and that contain no cyclic products. Using the anticommutativity of the generators, we may insist that $i_1 < i_2 < \cdots < i_k$. Using a rank argument, Jensen, McCammond, and Meier infer that this generating set is minimal. We conclude the following lemma. 

\begin{lem} \cite{JMM}. \label{CohomBasis}
 The cohomology group $H^k(\PS_n;\Z)$ has basis $$\{ \a_{i_1, j_1}^* \wedge \a_{i_2, j_2}^* \wedge \cdots \wedge \a_{i_k, j_k}^* \}$$ where $i_1 < i_2 < \cdots < i_k$, and where no permutations of the factors of $\a_{i_1, j_1}^* \wedge \a_{i_2, j_2}^*\wedge \cdots \wedge \a_{i_k, j_k}^* $ contains a cyclic product. 
\end{lem}

\begin{figure}[h] 
\begin{center}
 \includegraphics[scale=0.65]{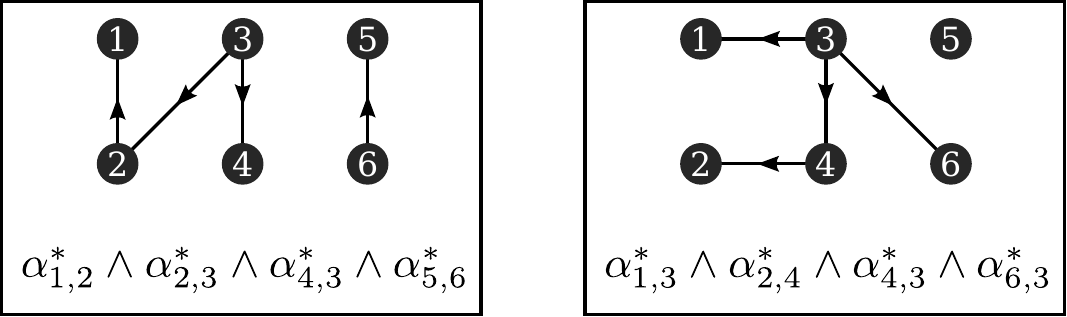}
 \caption{ {\small Two directed graphs and their corresponding basis elements in $H^4(\PS_6;\Q)$}}
 \label{fig:TwoForestExamples}
\end{center}
\end{figure}

Jensen, McCammond, and Meier count these basis elements by observing that they are indexed by directed graphs on $[n]$, where the graph contains an edge $i \leftarrow j$ whenever the basis element includes the factor $\a_{i,j}^*$. Two examples for $n=6$ and $k=4$ are given in \fullref{fig:TwoForestExamples}.

Jensen, McCammond, and Meier observe that, because each $i_m$ is distinct, each vertex is the target of at most one edge. This implies that graph cannot contain any cycles, as these would correspond to cyclic products. Thus, the connected components of the graph are trees, that is, the graph is a forest. Moreover, since each vertex is the target of at most one edge, each tree's digraph structure is determined by a choice of root. Since there are $k$ edges, the forest has $(n-k)$ connected components. In summary:

\begin{lem} \cite{JMM}. \label{TreeDecomp} 
 The cohomology group $H^k(\PS_n;\Z)$ is the free abelian group on the set of forests on $[n]$ of $(n-k)$ rooted trees.
\end{lem}

\noindent This description will be a key ingredient in our proofs of \fullref{WnStable} and \fullref{SnStable}.

\section{Representation Stability} \label{DefinitionRepStability}

The concept of representation stability was introduced by Church and Farb in \cite{RepStability}. Given a sequence of groups $\{G_n\}_{n \in \N}$  with inclusions $G_n \hookrightarrow G_{n+1}$, and a sequence $\{V_n \}_{n \in \N}$ of $G_n$--representations, we can describe a form of stability for $\{ V_n \}_{n \in \N}$ in terms of its decomposition into irreducible representations. Before giving a precise definition of representation stability, we introduce some notation used in \cite{RepStability}. \\ 

We will specialize to the cases that the group $G_n$ is either the symmetric group $S_n$, or the hyperoctahedral group $W_n$. The irreducible representations of $S_n$ over $\Q$ are classified by partitions $\y$ of $n$, where $\y = ( \y_0, \y_1, \ldots, \y_d)$ with $\y_0 \geq \y_1 \geq \cdots \geq \y_d$ and $\y_0 + \y_1 + \cdots + \y_d=n$. There is a natural identification of the irreducible representations of $S_n$ with certain irreducible representations of $S_{n+1}$, attained by associating the partition $\y = ( \y_0, \y_1, \ldots, \y_d)$ of $n$ with the partition $( \y_0 +1 , \y_1, \ldots, \y_d)$ of $(n+1)$. To emphasize this identification, from here on, we will describe a partition of $n$ by a partition $\y \vdash k$, $\y = (\y_1, \ldots, \y_d)$, with the first part $\y_0 = (n-k)$ left implicit. This notation is now independent of $n$ and gives a well-defined partition of $n$ whenever $n \geq (\y_1 + k)$. We call $$\y[n] := (n - k, \y_1, \ldots, \y_d)$$ the \emph{padded partition} of $\y$.  \\

The rational irreducible representations of the hyperoctahedral group $W_n$ are described in (for example) Geck--Pfeiffer \cite{GeckPfeiffer}. They are classified by \emph{double partitions} of $n$, that is, ordered pairs of partitions $\y = (\y^+, \y^-)$ with $|\y^+| + |\y^-| =  n$. These representations are constructed as follows: Given $\y \vdash n$ and $V(\y)$ the associated representation of $S_n$, we define $V( \y, (0))$ to be the $W_n$--representation pulled back from $V(\y)$ under the surjection $\pi \co W_n \twoheadrightarrow S_n$. Let $\varepsilon$ be the one-dimensional ``sign'' representation of $W_n$, where the canonical generators of $(\Z/2\Z)^n$ act by $(-1)$, and elements of $S_n$ act trivially. Define $$V( (0), \y)  : = V( \y, (0)) \otimes \varepsilon.$$ Next, we recall the following definition:
\begin{defn}
 Let $G$, $H$ be groups. Let $V$ be a rational representation of $G$, and $U$ a rational representation of $H$. Then the \emph{external tensor product} $V\boxtimes U$ denotes the $(G \times H)$--representation given by the action of $G \times H$ on the vector space $V \otimes U$ given by linearly extending the map $$(g,h) \cdot (v\otimes u) = ( g\cdot v ) \otimes ( h \cdot u ) \qquad \text{ for all } (g,h)\in G \times H \text{ and } v\otimes u \in V \otimes U.$$ 
\end{defn}

Recall additionally that, given a subgroup $H$ of a finite group $G$ and an $H$--representation $V$, the \emph{induced representation} $\Ind_H^G(V)$ is the $\Q G$--module $\Q G \otimes_{\Q H} V$. Equivalently, $\Ind_H^G(V)$ is the $G$--representation constructed as a direct sum $\Ind_H^G(V)=\bigoplus_{\s \in G/H} \s V$ of a translate of $V$ for each coset $\s H$. \\

 Then, generally, for $\y^+ \vdash k$ and $\y^- \vdash (n-k)$, we define $$ V(\y^+, \y^-) := \Ind_{W_k \times W_{n-k}}^{W_n} V(\y^+, (0)) \boxtimes V((0),\y^-). $$  Again, given a double partition $\y=(\y^+,\y^-)$ of $m$ with $|\y^-|=l$, we define the padded partition $$\y[n] := (\y^+[n-l], \y^-).$$ 

 For either $S_n$ or $W_n$, denote the irreducible representation associated to $\y[n]$ by $V(\y)_n$, and denote the multiplicity of $V(\y)_n$ in a representation $V_n$ by $c_{\y,n}(V_n)$. Using this notation, Church and Farb make the following definitions: 

\begin{defn}
A sequence $\{ V_n \}_{n \in \N}$ of $G_n$--representations with maps $\phi_n \co V_n \to V_{n+1}$ is called \emph{consistent} if $\phi_n$ is equivariant with respect to the $G_n$--action on $V_n$ and the $G_n$--action on $V_{n+1}$ induced by the natural inclusion $G_n \hookrightarrow G_{n+1}$. \end{defn}

\begin{defn}
 A consistent sequence $\{V_n\}_{n \in \N}$ of $G_n$--representations with maps $\phi_n \co V_n \to V_{n+1}$ is called \emph{representation stable} if it satisfies the following properties: \\

I. \textbf{\emph{Injectivity.}} The maps $\phi_n \co V_n \to V_{n+1}$ are injective, for $n$ sufficiently large,.\\

II. \textbf{\emph{Surjectivity.}} The space $V_{n+1}$ is the $G_{n+1}$--span of $\phi_n(V_n)$, for $n$ sufficiently large. \\

III. \textbf{\emph{Multiplicities.}} In the decomposition $$V_n = \bigoplus_{\y}c_{\y,n}V(\y)_n$$ of $V_n$ into irreducible $G_n$--representations, for each $\y$, the coefficient $c_{\y,n}$ is eventually independent of $n$. 
\end{defn}

\begin{defn}
 A representation stable sequence $\{ V_n \}_{n \in \N}$ is \emph{uniformly} representation stable if the multiplicities $c_{\y,n}(V_n)$ become constant at some $N \in \N$ not depending on $\y$. 
\end{defn}

Our goal is now to prove that the cohomology groups $H^*(\PS_n, \Q)$ are uniformly representation stable as $W_n$--representations.

\section{The Cohomology Groups of the Pure String Motion Group are Uniformly Representation Stable} \label{MainThm}

The action of $W_n$ on $\PS_n$ by conjugation induces an action on the cohomology: 

\[
 \p_i \cdot \a_{j,k}^* = \left\{ \begin{array}{lll}
        - \a_{j,k}^* & (k=i) &  \\  
          \a_{j,k}^*  & (k \neq i) &  \end{array} \right. \qquad 
\] \[
\t_i \cdot \a_{j,k}^* =  \a_{ \sigma(j), \sigma(k) }^* \quad \text{where } \sigma = (i \; \; \; i+1) \in S_n.
\]
With this action, we will consider the rational cohomology group $H^k(\PS_n; \Q)$ as a $W_n$--representation. \\

The maps $\PS_{n+1} \to \PS_n$ obtained by `forgetting' the circle $\g_{n+1}$ induce $W_n$--equivariant maps  
\begin{align*}
\phi_n \co H^*(\PS_n ; \Q) &\longrightarrow H^*(\PS_{n+1} ; \Q)\\
\a_{i,j}^* & \longmapsto \a_{i,j}^* 
\end{align*}
Thus, $\{(H^k(\PS_n ; \Q), \phi_n) \}_{n \in \N}$ is a consistent sequence. \\

We will use the strategy of Church and Farb for the proof of Theorem 4.1 in \cite{RepStability}, the uniform representation stability of the cohomology groups of the pure braid group, to prove the following:

\begin{thm}\label{WnStable}
For each fixed $k \geq 0$, the sequence of $W_n$--representations $$\{H^k(\PS_n ; \Q)\}_{n \in \N}$$ is uniformly representation stable with respect to the maps $$\phi_n \co H^k(\PS_n ; \Q) \to H^k(\PS_{n+1} ; \Q)$$ induced by the `forgetful' map $\PS_{n+1} \to \PS_n$. The sequence stabilizes once $n\geq 4k$.
\end{thm}

Injectivity follows easily from the definition of the maps $\phi_n$, since a basis element $$\a_{i_1, j_1}^* \wedge \cdots \wedge \a_{i_k, j_k}^* \in H^k(\PS_n ; \Q)$$ maps to the basis element $$\a_{i_1, j_1}^* \wedge \cdots \wedge \a_{i_k, j_k}^* \in H^k(\PS_{n+1} ; \Q).$$ For the surjectivity criterion, observe that for $n \geq 2k$ and any basis element $v= \a_{i_1, j_1}^* \wedge \cdots \wedge \a_{i_k, j_k}^*$ of $H^k(\PS_n ; \Q)$, there must exist some number in $[n+1]$ that does not appear among the indices of $v$. Then, for some suitable choice of permutation $\sigma \in S_{n+1} \subseteq W_{n+1}$, the basis element $\sigma \cdot v$ does not have $(n+1)$ as an index, and so is contained in the image $\phi_n(H^k(\PS_n ; \Q))$. \\

To prove the multiplicities criterion, we will use a combinatorial result concerning certain induced representations. In \cite{Hemmer}, Hemmer proved the following for the symmetric group: 

\begin{lem} \label{IndSn} \cite[Theorem 2.4]{Hemmer}.
Fix $r  \geq 1$. Let $H$ be any subgroup of $S_{r }$ and $V$ an $H$--representation. Let $\Q$ denote the trivial representation of $S_{n-r }$. Then the sequence of $S_n$--representations $\{ \Ind^{S_n}_{H \times S_{n-r }}(V \boxtimes \Q) \}_{n \in \N}$ is uniformly representation stable, stabilizing for $n \geq 2r $. 
\end{lem} 

In \cite[Theorem 4.6]{RepStability}, Church and Farb adapt the proof to the hyperoctahedral group. We reproduce their proof, with more explicit attention paid to the bounds on the stable range.

\begin{lem} \label{IndWn} \cite{RepStability}.
Fix $r  \geq 1$. Let $H$ be any subgroup of $W_{r }$, and $V$ an $H$--representation. Let $\Q$ denote the trivial representation of $W_{n-r }$. Then the sequence of $W_n$--representations $\{ \Ind^{W_n}_{H \times W_{n-r }}(V \boxtimes \Q) \}_{n \in \N}$ is uniformly representation stable, stabilizing for $n \geq 2r$.
\end{lem} 

\begin{proof}[Proof of \fullref{IndWn}]

We first observe that $$ \Ind^{W_n}_{H \times W_{n-r }}(V \boxtimes \Q) = \Ind^{W_n}_{W_r \times W_{n-r}}\big( (\Ind^{W_r}_H V) \boxtimes \Q \big),$$ and so it would suffice to prove the theorem for an irreducible representation $V(\y^+, \y^-)$ of $W_{r}$, and then apply the result to each irreducible factor of $\Ind^{W_r}_H V$. Throughout this proof, for notational simplicity, we refer to irreducible representations by the full padded partitions.\\

By Lemma 6.1.3 of Geck--Pfeiffer \cite{GeckPfeiffer}, 
$$\Ind^{W_n}_{W_r \times W_{n-r}} \big( V(\y^+, \y^-) \boxtimes V(\mu^+, \mu^-) \big) = \bigoplus_{(\nu^+, \nu^-)} C_{\y^+, \mu^+}^{\nu^+} C_{\y^-, \mu^-}^{\nu^-} V(\nu^+, \nu^-) $$
where $C_{\y, \mu}^{\nu}$ denotes the Littlewood-Richardson coefficient. According to the classification of irreducible hyperoctahedral representations, the trivial representation of $W_{n-r}$ is $V((n-r),(0))$, and so 
$$\Ind^{W_n}_{W_r \times W_{n-r}} \big( V(\y^+, \y^-) \boxtimes \Q \big) = \bigoplus_{(\nu^+, \nu^-)} C_{\y^+, (n-r)}^{\nu^+} C_{\y^-, (0)}^{\nu^-} V(\nu^+, \nu^-). $$

It is a standard result that 
\[
 C_{\y^-, (0)}^{\nu^-} = \left\{ \begin{array}{ll}
        1 \quad \text{if } \nu^- = \y^- & \\
        0 \quad \text{if } \nu^- \neq \y^-. & \end{array}  \right.  \] 

(See, for example, Fulton \cite{FultonYoungTableaux} for the combinatorics of the Littlewood-Richardson coefficients.) \\

Then, 
\begin{align*}
 \Ind^{W_n}_{W_r \times W_{n-r}} \big( V(\y^+, \y^-) \boxtimes \Q \big) &= \bigoplus_{\nu^+} C_{\y^+, (n-r)}^{\nu^+} V(\nu^+, \y^-) \\
 & = \bigoplus_{\nu^+} V(\nu^+, \y^-)
\end{align*}
 where the final sum is taken over all Young diagrams $\nu^+$ that are constructed by adding $(n-r)$ boxes to the Young diagram $\y^+$, with no two boxes added to the same column. When $n > 2r$, then $n-r > r$, so $(n-r)$ is necessarily greater than the number of columns of $\y^+$, which implies that at least one box must be added to the first row. Thus, when $n > 2r$, the function `removing a box from the first row' gives a bijection between the admissible partitions of $n$ and those of $(n-1)$, a bijection which identifies a double partition $\nu[n]$ with $\nu[n-1]$. It follows that the decomposition stabilizes once $n \geq 2r$. \end{proof}

Having proved \fullref{IndWn}, our objective now is to decompose the representations $H^k(\PS_n ; \Q)$ into a sum of subrepresentations of the form $\Ind^{W_n}_{H \times W_{n-r }}(V \boxtimes \Q )$. 

\begin{proof}[Proof of \fullref{WnStable}] 
 By \fullref{TreeDecomp}, $H^k(\PS_n ; \Q)$ is the rational vector space over a basis indexed by the forests on $[n]$ with $(n-k)$ connected components. Let $\cP$ be a partition of $[n]$ into $(n-k)$ subsets; we call its subsets \emph{blocks}. Define $H^{\cP}(\PS_n; \Q) \subseteq H^k(\PS_n ; \Q)$ to be the $\Q$--span of the set of those forests whose connected components are the blocks of $\cP$. Define $\bcP$ to be the partition of $n$ associated to $\cP$, given by the block sizes.\\

Consider, for example, the case where $n=6$, $k=2$, and $\cP = \{ \{1\}, \{2,4\}, \{3,5\}, \{6\} \}$. Then $H^{\cP}(\PS_n; \Q)$ is the four-dimensional subspace spanned by the basis elements shown in \fullref{fig:OnePartitionExample}, and $\bcP = (2,2,1,1).$ \\ 

\begin{figure}[h] 
\begin{center}
 \includegraphics[scale=0.5]{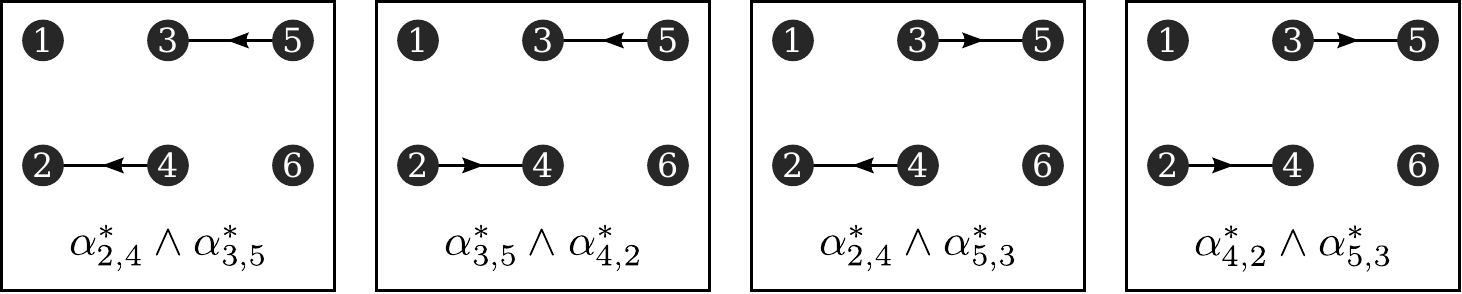}
 \caption{ {\small A basis for $H^{\{ \{1\}, \{2,4\}, \{3,5\}, \{6\} \}}(\PS_6; \Q)$} }
\label{fig:OnePartitionExample}
\end{center}
\end{figure}

Given any $n$ and any partition $\cP$ of $[n]$, we can compute the action of $W_n$ on $H^{\cP}(\PS_n; \Q)$. An element $\p_i \in (\Z/2\Z)^n \subseteq W_n$ either fixes or negates each basis element of $H^k(\PS_n ; \Q)$, and so stabilizes $H^{\cP}(\PS_n; \Q)$. An element $\t_i \in S_n \subseteq W_n$ permutes the labels on the vertices, and so the $S_n$--orbit of $H^{\cP}(\PS_n; \Q)$ are all subspaces $H^{\cP'}(\PS_n; \Q)$ of $H^k(\PS_n ; \Q)$ with $\overline{\cP'} = \bcP$. We therefore have a decomposition: 
\begin{align*}
 H^k(\PS_n ; \Q) &= \bigoplus_{  \substack{ \text{partitions } \cP \text{ of } [n] \\  \text{ with } (n-k) \text{ blocks} }  } H^{\cP}(\PS_n) \\ 
&= \bigoplus_{  \substack{  \mu \vdash n \\  |\mu|=(n-k) }  } \; \bigoplus_{ \{ \cP \; | \; \bcP = \mu \} }  H^{\cP}(\PS_n) \\ 
&= \bigoplus_{ \substack{  \mu \vdash n \\  |\mu|=(n-k) } } \Ind_{\Stab(\cP_{\mu})}^{W_n} H^{\cP_{\mu}}(\PS_n)
\end{align*}

 for some fixed $\cP_{\mu}$ with $\overline{\cP_{\mu}} = \mu$. \\ 

For any $\mu$, the stabilizer $\Stab(\cP_{\mu}) \subseteq W_n$ is generated by all elements $\p_i$, by the permutations that permute the numbers within a block of $\cP_{\mu}$, and by the permutations that interchange two blocks of $\cP_{\mu}$ of the same size. If $\cP_{\mu}$ has $(n-r)$ blocks of size one, then we may decompose the stabilizer as $\Stab(\cP_{\mu}) \cong H \times W_{n-r}$,  where $H$ is a subgroup of $W_r$, and $W_{n-r}$ acts by permuting the singleton blocks. Since the singleton blocks of $\cP_{\mu}$ are those numbers that do not appear as an index in any basis element of $H^{\cP_{\mu}}(\PS_n; \Q)$, this direct factor $W_{n-r}$ acts trivially on $H^{\cP_{\mu}}(\PS_n; \Q)$. \\

 We note moreover that, whenever $k > 0$, not all vertices can be singletons, and so we satisfy the hypothesis of \fullref{IndWn} that $r \geq 1$.\\

The maps $$\phi_{n-1} \co H^*(\PS_{n-1} ; \Q) \to H^*(\PS_{n} ; \Q)$$ map $ \a_{i_1, j_1}^* \wedge \cdots \wedge \a_{i_k, j_k}^* \in H^k(\PS_{n-1} ; \Q)$ to $ \a_{i_1, j_1}^* \wedge \cdots \wedge \a_{i_k, j_k}^* \in H^k(\PS_{n} ; \Q)$, and accordingly, they map a rooted forest of $((n-1)-k)$ trees on $[n-1]$ to the rooted forest of $(n-k)$ trees on $[n]$ obtained by adding the singleton vertex $n$. For each partition $\cP$ of $[n-1]$, the map $\phi_{n-1}$ induces an isomorphism of $W_{n-1}$--representations
 \begin{equation*} H^{\cP}(\PS_{n-1}; \Q ) \cong  H^{\cP \cup \{ \{ n \} \} }(\PS_{n}; \Q ). \end{equation*} 
When $n > 2k$, any partition $\cP^+$ of $[n]$ into $(n-k)$ blocks must have at least one singleton block. Thus $H^{\cP^+}(\PS_{n}; \Q )$ is in the $W_n$--orbit of $(\cP \cup \{ \{ n \} \})$ for some partition $\cP$ of $[n-1]$. It follows that:
\begin{align*}
 H^k(\PS_{n} ; \Q) &= \bigoplus_{ \substack{  \mu \vdash (n) \\  |\mu|=(n-k) } } \Ind_{\Stab(\cP_{\mu})}^{W_{n}} H^{\cP_{\mu}}(\PS_{n}) \\
&= \bigoplus_{ \substack{  \nu \vdash n-1 \\  |\nu|=(n-1-k) } } \Ind_{\Stab(\cP_{\nu} \cup \{ \{ n \} \} )}^{W_{n-1}} H^{\cP_{\nu} \cup \{ \{ n \} \} }(\PS_{n}) 
\end{align*}
for some choices of $\cP_{\mu}$ and $\cP_{\nu}$ with $\overline{\cP_{\mu}} = \mu$, $\overline{\cP_{\nu}} = \nu$. The number of direct summands in this decomposition therefore stabilizes once $n \geq 2k$.\\

For given $\nu$, we may choose $\cP_{\nu}$ so that its $((n-1)-r)$ singleton blocks are $\{ \{r+1\}, \{r+2\}, \ldots, \{n-1\} \}$, and denote the complement, a partition of $[r]$, by $\cP^c$. Since $\Stab(\cP_{\nu}\cup \{ \{ n \} \} ) \cong H \times W_{n-r}$ as above, we can iteratively write:
\begin{align*} \Ind_{\Stab(\cP_{\nu}  \cup \{ \{ n \} \} )}^{W_n} H^{\cP_{\nu}  \cup \{ \{ n \} \} }(\PS_{n}) &= \Ind_{H \times W_{n-r} }^{W_n} H^{\cP_{\nu}  \cup \{ \{ n \} \} }(\PS_{n}) \boxtimes \Q \\
 &= \Ind_{H \times W_{n-r}}^{W_n} H^{\cP_{\nu} } (\PS_{n-1}) \boxtimes \Q \\  
 &= \ldots \\ 
 &= \Ind_{H \times W_{n-r}}^{W_n} H^{\cP^{c} }(\PS_{r}) \boxtimes \Q. \end{align*}

We may conclude that the sequence $\{H^k(\PS_{n} ; \Q)\}_{n \geq 2k}$ of $W_n$--representations decomposes into a finite sum of sequences of the form described in \fullref{IndWn}. \\

Among all partitions of $[n]$ into $(n-k)$ blocks, the number of singleton blocks is minimal in the partition with $k$ blocks of two elements, and $(n-2k)$ singleton blocks. This implies that the maximal value of $r$ is $2k$, and so we can ensure that $n \geq 2r$ whenever $n \geq 4k$.  \fullref{WnStable} now follows from \fullref{IndWn}. \end{proof}

We notice that we may use this same method of proof to prove representation stability when we restrict the action of $W_n$ to its subgroup $S_n$.

\begin{thm}\label{SnStable}
For each fixed $k \geq 0$, sequence of groups $$\{H^k(\PS_n ; \Q)\}_{n \in \N},$$ considered as $S_n$--representations, is uniformly representation stable with respect to the maps $$\phi_n \co H^k(\PS_n ; \Q) \to H^k(\PS_{n+1} ; \Q)$$ induced by the `forgetful' map $\PS_{n+1} \to \PS_n$. The sequence stabilizes once $n \geq 4k$.
\end{thm}
\begin{proof}[Proof of \fullref{SnStable}] Since this proof is so similar to that of \fullref{WnStable}, we merely give a sketch. We may use the same decomposition of the cohomology groups as above, noting that 
\begin{align*}
 H^k(\PS_n ; \Q) &= \bigoplus_{ \substack{  \mu \vdash n \\  |\mu|=(n-k) } } \Ind_{\Stab(\cP_{\mu})}^{S_n} H^{\cP_{\mu}}(\PS_n)
\end{align*}
 for some $\cP_{\mu}$ with $\overline{\cP_{\mu}} = \mu$. Here, if $\cP_{\mu}$ has $m_j$ blocks of size $j$, then its stabilizer $\Stab(\cP_{\mu}) \subseteq S_n$ decomposes as
\begin{align*}
\Stab(\cP_{\mu}) = \bigoplus_{j=1}^n  S_j \wr S_{m_j} 
\end{align*}
 where $S_j$ permutes the numbers within a block of size $j$, and $S_{m_j}$ permutes the $m_j$ blocks. In particular, if $\cP_{\mu}$ has $(n-r)$ blocks of size one, then we can write $\Stab(\cP_{\mu}) \cong H \times S_{n-r}$, where $S_{n-r}$ permutes the singleton blocks. Then, using \fullref{IndSn} as above, we may conclude the theorem. 
\end{proof}

These results raise the question of whether it is true in general that a uniformly representation stable sequence of $W_n$--representations will restrict to a uniformly representation stable sequence of $S_n$--representations. This is indeed the case, and will be shown in a forthcoming paper. However, \fullref{SnStable} gives a better stable range than the general case.

\section{The Cohomology of the String Motion Group} \label{StringMotionGroup}

Consider a short exact sequence of groups 

\begin{center}

$\xymatrix{ & 1 \ar[r] & N  \ar[r] & G  \ar[r] & Q  \ar[r] & 1 }$

\end{center}

\noindent with finite quotient $Q$. Recall that, by a transfer argument, the restriction map $H^k(G; \Q) \to H^k(N; \Q)$ identifies the cohomology of $G$ with the ring of invariants $H^k(N; \Q)^Q$ under the natural action of $Q$. (See, for example, Brown \cite[Chapter 3.9]{Brown}.) It therefore follows from the exact sequence

\begin{center}

$\xymatrix{ & 1 \ar[r] & \PS_n  \ar[r] & \S_n  \ar[r] & W_n  \ar[r] & 1 }$

\end{center}

\noindent that the cohomology of $H^k(\S_n; \Q)$ is isomorphic to the isotypic component of trivial $W_n$--representations $V((0), (0))_n$ in $H^k(\PS_n ; \Q)$. \\

 The representation stability of $H^k(\PS_n; \Q)$ implies, in particular, that the multiplicity of the trivial representation $V((0), (0))_n$ stabilizes for $n$ sufficiently large, and therefore that the dimension of $H^k(\S_n; \Q)$ is eventually constant. In fact, even the integral homology groups of $\S_n$ are eventually independent of $n$, a result proved by Hatcher and Wahl \cite[Corollary 1.2]{HatcherWahl}.\\ 

Hatcher and Wahl realize $\S_n$ as follows: Let  $P := S^1 \times D^2$ be the product of a $1$--sphere $S^1$ and the $2$--disk $D^2$. Let $N_n := D^3 \# P \# \cdots \# P$ be the connected sum of the $3$--disk $D^3$ with $n$ copies of $P$. Then $\S_n$ is the quotient of the mapping class group $\Gamma(N_n, \partial D^3) :=\pi_0 \Diff(N_n \; \rel \; \partial D^3)$ by the subgroup generated twists along spheres and disks with boundary in $\partial N_n - \partial D^3$. Homological stability for $\S_n$ is a consequence of their stability theorems for mapping class groups of $3$--manifolds.\\

From our description of $H^k(\PS_n; \Q)$, we can say more: for $k \geq 1$, the trivial representation does not occur, and so the rational cohomology groups of $\S_n$ vanish.  

\begin{thm} \label{SigmaVanish}
 For $k \geq 1$, the rational cohomology groups $H^k(\S_n; \Q)=0$. In particular, the integral homology and cohomology groups of $\S_n$ are finite. 
\end{thm}

\begin{proof}[Proof of \fullref{SigmaVanish}]
 For a $G$--representation $V$, the averaging map 
\begin{align*}
\varphi \co  V &\to V \\
       u & \mapsto \frac{1}{|G|} \sum_{g \in G} g \cdot u
\end{align*}
gives a linear projection of $V$ onto the subspace fixed by $G$. Thus, to prove that $H^k(\PS_n; \Q)$ contains no copies of the trivial representation, it suffices to show that $\varphi(v)=0$ for every basis element $v$. \\

Suppose that, for some fixed basis element $v$, there were an involution $\omega \in W_n$ such that $\omega \cdot v = -v$. Then by partitioning the sum by left cosets of $H=\{\omega, 0\}$, we find that $\varphi(v)$ vanishes:
\begin{align*}
\varphi(v) = \frac{1}{|G|} \sum_{g \in W_n} g \cdot v =  \frac{1}{|G|} \sum_{gH \subseteq W_n} \big( g \cdot v + g \omega \cdot v \big) = \frac{1}{|G|} \sum_{gH \subseteq W_n} \big( g \cdot (v -v) \big) =0. 
\end{align*}

For any basis element $v = \a_{i_1, j_1}^* \wedge \a_{i_2, j_2}^* \wedge \cdots \wedge \a_{i_k, j_k}^*$ of $H^k(\PS_n; \Q)$, we will find such an involution. Suppose that there is some number $j \in [n]$ that appears as a second index for an odd number of exterior factors of $v$. Then $\rho_j \cdot v = -v$, which implies that $\varphi(v)=0$. Next, suppose that there is no such number $j$. This means that, if we take the rooted forest associated to $v$ as described in \fullref{CohomologyPSn}, every vertex has even outdegree.  Choosing any tree in the forest with at least one edge, we can find two terminal vertices $p$ and $q$ that have the same direct predecessor $r$. This implies $v$ includes the factors $\a_{p, r}^*$ and $\a_{q, r}^*$, but that the numbers $p$ and $q$ do not appear as any other index. Thus, the transposition $ \omega = (p \; q) \in S_n$ induces an odd permutation on the factors of the exterior product $v$, and $\varphi(v)=0$. This completes the proof that $H^k(\S_n; \Q)=0$ for $k \geq 1$. \\

 The work of Collins \cite{Collins} (later Jensen and Wahl \cite{JensenWahl}) implies that  $\PS_n$ is of type $FP_{\infty}$. Since $[\S_n : \PS_n] = 2^n n! < \infty$, it follows that $\S_n$ is of type $FP_{\infty}$; a standard reference is Brown \cite[Chapter 8]{Brown}. The integral homology and cohomology groups are finitely generated and torsion; they are finite groups. 
\end{proof}

An approach to computing these groups is developed in Griffin \cite{Griffin}.

\section{The Cohomology of the Braid-Permutation Group} \label{BraidPermutation}

This section originated in conversations with Allen Hatcher. It also owes to insights from Tom Church. \\ 

Consider the group $\S_n^+ = \langle \tau_i, \alpha_{i,j} \rangle$ of string motions that preserve the orientations of the circles. This group has been called the \emph{braid-permutation group}, since it is generated by the classical braid group $B_n= \langle \a_{i,i+1} \t_i \rangle$ and the symmetric group $S_n=\langle \tau_i \rangle$. Fenn, Rimanyi, and Rourke \cite{FR, FRR} describe the group as a generalized braid group in which some crossings are ``welded''. They give a presentation for the group, and also show that it is isomorphic to the automorphism group of the so-called \emph{free quandle}. In \cite{Vershinin}, Vershinin studies the algebraic topology of the braid-permutation group. Applications of the group to mathematical physics are described by Baez, Crans, and Wise in \cite{Exotic}. Brendle and Hatcher \cite{BrendleHatcher} give a geometric description of $\S_n^+$ as the \emph{untwisted ring group}, the fundamental group of the $4n$--manifold of configurations of $n$ Euclidean circles in $\
R^3$ lying parallel to a fixed plane.\\

 The symmetries of the groups $H^k(\PS_n; \Q)$ encode the structure of the cohomology of  $\S_n^+$. From the short exact sequence 

\begin{center}

$\xymatrix{ & 1 \ar[r] & \PS_n  \ar[r] & \S_n^+  \ar[r] & S_n  \ar[r] & 1 }$

\end{center}

\noindent we see that the cohomology group $H^k(\S_n^+; \Q)$ is the group of invariants $H^k(\PS_n; \Q)^{S_n}$. The representation stability of the trivial $S_n$--representation $V(0)_n$ in \fullref{SnStable} therefore implies classical rational stability for $H^k(\S_n^+; \Q)$.

\begin{cor} \label{BPStable}
 The braid-permutation group $\S_n^+$ of orientation-preserving string motions is rationally cohomologically stable in the classical sense. For each $k \geq 0$, $H^k(\S_n^+; \Q) \cong H^k(\S_{n+1}^+; \Q)$ once $n\geq 4k$.
\end{cor}

\begin{rem}
 In contrast to the case with the $\PS_n$, we cannot construct a map $\S_{n+1}^+ \to \S_n^+$ by ``forgetting'' a circle, as there is no consistent way to chose among $(n+1)$ unlabeled circles. Since $\R^3$ is open, however, there are maps $\S_{n}^+ \to \S_{n+1}^+$ that may be constructed by `introducing' a circle far from the origin of $\R^3$. These maps induce the isomorphisms in \fullref{BPStable}.
\end{rem}

For low degree $k$, we can use the symmetries of $H^k(\PS_n; \Q)$ to directly compute the $S_n$--invariant subspace: \\ 

For $n \geq 2$, the first cohomology $H^1(\S_n^+; \Q) = \Q$ is one-dimensional, corresponding to the trivial representation spanned by
\begin{align*} \sum_{i \neq j \in [n]} \a_{i,j}^* \in H^1(\PS_n; \Q).  \end{align*} 

For $n \geq 3$, the second cohomology $H^2(\S_n^+; \Q) = \Q$ is also one-dimensional, corresponding to the trivial representation spanned by
\begin{align*} \sum_{i, j, k \in [n] \text{ distinct}} \a_{i,j}^* \wedge \a_{j,k}^* \in H^2(\PS_n; \Q).  \end{align*}  

For $n \geq 5$, the cohomology $H^3(\S_n^+; \Q) = \Q^3$ corresponding to the trivial representations spanned by
\begin{align*} \sum_{i, j, k,l \in [n] \text{ distinct}} \a_{i,j}^* \wedge \a_{j,k}^* \wedge \a_{k,l}^* & \in H^3(\PS_n; \Q),\\
 \sum_{i, j, k,l \in [n] \text{ distinct}} \a_{i,j}^* \wedge \a_{j,k}^* \wedge \a_{l,k}^* & \in H^3(\PS_n; \Q), \text{ and} \\ 
\sum_{i, j, k,l,m \in [n] \text{ distinct}} \a_{i,j}^* \wedge \a_{j,k}^* \wedge \a_{l,m}^* & \in H^3(\PS_n; \Q). \end{align*}

\begin{figure}[h]
\begin{center}
 \includegraphics[scale=0.64]{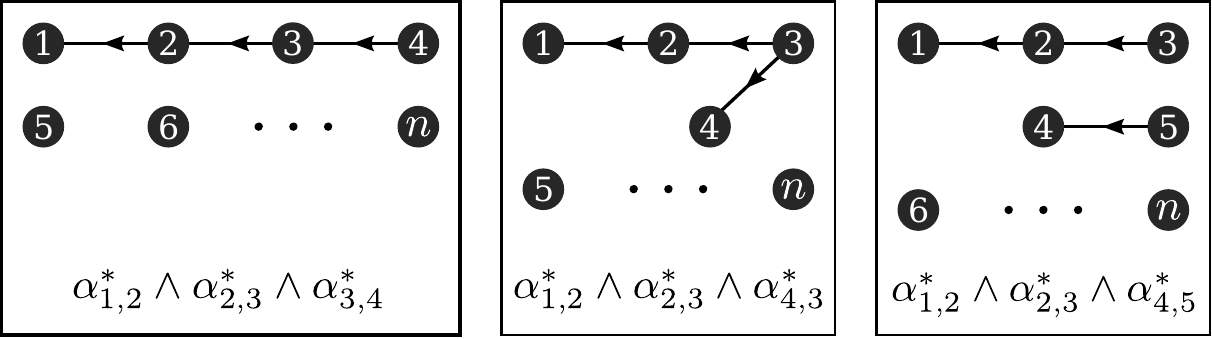}
\caption{ {\small Elements of $H^3(\PS_n; \Q)$ corresponding to a basis for $H^3(\S_n^+; \Q)$}}
 \label{fig:SigmaPlusH3Basis}
\end{center}
\end{figure}

These three vectors are obtained by summing over the $S_n$--orbits of each of the three basis elements shown (respectively) in \fullref{fig:SigmaPlusH3Basis}.

In general, each of these trivial reprsentations is spanned by the symmetrization of some basis element $v \in H^k(\PS_n; \Q)$. To ensure that this symmetrization does not vanish, we need $v \in H^k(\PS_n; \Q)$ with the property that no element of $S_n$ negates $v$. Taking the basis for $H^k(\PS_n; \Q)$ described in \fullref{CohomBasis} and \fullref{TreeDecomp}, the symmetries of a basis element may be seen from the associated rooted forests. Specifically, the trivial representations will correspond to forests such that any permutation of the vertices preserving the labelled, rooted graph structure will induce an even permutation of the factors of the associated basis element. Such forests exist for all $1 \leq k \leq (n-1)$; consider, for example, the simple path $ 1 \leftarrow 2 \leftarrow \cdots \leftarrow k \leftarrow (k+1)$ along with isolated vertices $(k+2)$ through $n$. Thus, for all $n$, $H^k(\S_n^+; \Q)$ is nontrivial in all dimensions $0 \leq k \leq (n-1)$. \\

In \cite[Example 6.2]{Griffin}, Griffin gives a method of computing the integral homology of $\S_n^+$ from the structure of these rooted forests. \\ 

In \cite{Vershinin}, Vershinin proves that the plus construction of the classifying space of the infinite braid-permution group is homotopic to a component of $\Omega^{\infty}S^{\infty} \times S^1 \times Y$, where $Y$ is an infinite loop space. A better understanding of the combinatorics of the rooted forests could illuminate the cohomological structure of this space. \\

\bibliographystyle{gtart}
\bibliography{MasterBibliography}

\begin{thebibliography}{}
\providecommand\bibmarginpar{\leavevmode\marginpar}
\def\urlstyle#1{{\tt #1}}

\bibitem{Exotic}
\textbf{J\,C Baez}, \textbf{D\,K Wise}, \textbf{A\,S Crans}, \emph{Exotic
  statistics for strings in 4d BF theory}, Advances in Theoretical and
  Mathematical Physics 11 (2007) 707--749

\bibitem{DualityGroup}
\textbf{N Brady}, \textbf{J McCammond}, \textbf{J Meier}, \textbf{A Miller},
  \emph{The Pure Symmetric Automorphisms of a Free Group Form a Duality Group*
  1}, Journal of Algebra 246 (2001) 881--896

\bibitem{BrendleHatcher}
\textbf{T Brendle}, \textbf{A Hatcher}, \emph{Configuration spaces of rings and
  wickets}, Commentarii Mathematici Helvetici  (2011)

\bibitem{Brown}
\textbf{K\,S Brown}, \emph{Cohomology of groups}, Springer (1982)

\bibitem{BrownsteinLee}
\textbf{A Brownstein}, \textbf{R Lee}, \emph{Cohomology of the group of motions
  of n strings in 3-space}, from: ``Mapping class groups and moduli spaces of
  Riemann surfaces: proceedings of workshops held June 24-28, 1991, in
  G{\"o}ttingen, Germany, and August 6-10, 1991, in Seattle, Washington...'',
  volume 160, Amer Mathematical Society (1993)  51

\bibitem{RepStability}
\textbf{T Church}, \textbf{B Farb}, \emph{{Representation theory and
  homological stability}}, Arxiv preprint arXiv:1008.1368  (2010)

\bibitem{Collins}
\textbf{D\,J Collins}, \emph{Cohomological dimension and symmetric
  automorphisms of a free group}, Commentarii Mathematici Helvetici 64 (1989)
  44--61

\bibitem{Dahm}
\textbf{D\,M Dahm}, \emph{A generalization of braid theory.}, Princeton Univ.
  Ph. D. thesis  (1962)

\bibitem{FR}
\textbf{R Fenn}, \textbf{R Rim{\'a}nyi}, \emph{Some remarks on the
  braid-permutation group}, Topics in Knot Theory 399 (1993) 57--68

\bibitem{FRR}
\textbf{R Fenn}, \textbf{R Rim{\'a}nyi}, \textbf{C Rourke}, \emph{The
  braid-permutation group}, Topology 36 (1997) 123--135

\bibitem{FultonYoungTableaux}
\textbf{W Fulton}, \emph{Young tableaux}, Cambridge Univ. Press (1997)

\bibitem{GeckPfeiffer}
\textbf{M Geck}, \textbf{G Pfeiffer}, \emph{Characters of finite Coxeter groups
  and Iwahori-Hecke algebras}, Oxford University Press, USA (2000)

\bibitem{Goldsmith}
\textbf{D\,L Goldsmith}, \emph{The theory of motion groups.}, The Michigan
  Mathematical Journal 28 (1981) 3--17

\bibitem{Griffin}
\textbf{J Griffin}, \emph{Diagonal complexes and the integral homology of the
  automorphism group of a free product}, Arxiv preprint arXiv:1011.6038  (2010)

\bibitem{HatcherWahl}
\textbf{A Hatcher}, \textbf{N Wahl}, \emph{Stabilization for mapping class
  groups of 3-manifolds}, Duke Mathematical Journal 155 (2010) 205--269

\bibitem{Hemmer}
\textbf{D\,J Hemmer}, \emph{{Stable decompositions for some symmetric group
  characters arising in braid group cohomology}}, Journal of Combinatorial
  Theory, Series A  (2010)

\bibitem{JMM}
\textbf{C\,A Jensen}, \textbf{J McCammond}, \textbf{J Meier}, \emph{The
  integral cohomology of the group of loops}, Geometry \& Topology 10 (2006)
  759--784

\bibitem{JensenWahl}
\textbf{C\,A Jensen}, \textbf{N Wahl}, \emph{Automorphisms of free groups with
  boundaries}, Algebr. Geom. Topol 4 (2004) 543--569

\bibitem{McCool}
\textbf{J McCool}, \emph{On basis-conjugating automorphisms of free groups},
  Canadian J. Math 38 (1986) 1525--1529

\bibitem{PettetFiniteness}
\textbf{A Pettet}, \emph{Finiteness properties for a subgroup of the pure
  symmetric automorphism group}, Comptes Rendus Mathematique 348 (2010)
  127--130

\bibitem{Rubinsztein}
\textbf{R\,L Rubinsztein}, \emph{On the group of motions of oriented, unlinked
  and unknotted circles in $\mathbb{R}^3$, I}, Preprint  (2002)

\bibitem{Vershinin}
\textbf{V\,V Vershinin}, \emph{On homological properties of singular braids},
  Transactions-American Mathematical Society 350 (1998) 2431--2456

\bibitem{Wattenberg}
\textbf{F Wattenberg}, \emph{Differentiable motions of unknotted, unlinked
  circles in 3-space}, Math. Scand 30 (1972) 107--135

\end{thebibliography}



\end{document}